\documentclass[11pt]{amsart}

\usepackage{amssymb,bbm,amscd}
\usepackage{graphicx}
\usepackage{a4wide}

\newcommand{\dd}{\,\mathrm{d}}
\newcommand{\ts}{\hspace{0.5pt}}

\begin{document}

\title[Diffraction of the Thue-Morse chain]{The singular continuous 
diffraction measure\\[2mm] of the Thue-Morse chain}

\author{Michael Baake}
\address{Fakult\"{a}t f\"{u}r Mathematik, Universit\"{a}t Bielefeld,\newline
\hspace*{\parindent}Postfach 100131, 33501 Bielefeld, Germany}
\email{mbaake@math.uni-bielefeld.de}

\author{Uwe Grimm}
\address{Department of Mathematics and Statistics,
The Open University,\newline 
\hspace*{\parindent}Walton Hall, Milton Keynes MK7 6AA, United Kingdom}
\email{u.g.grimm@open.ac.uk}

\begin{abstract}
  The paradigm for singular continuous spectra in symbolic dynamics
  and in mathematical diffraction is provided by the Thue-Morse chain,
  in its realisation as a binary sequence with values in $\{\pm
  1\}$. We revisit this example and derive a functional equation
  together with an explicit form of the corresponding singular
  continuous diffraction measure, which is related to the known
  representation as a Riesz product.
\end{abstract}

\maketitle
\thispagestyle{empty}

\noindent
The Thue-Morse chain can be defined via the primitive substitution rule
\[
    \varrho : \;
    \begin{array}{r@{\;}c@{\;}l}
    1 & \mapsto & 1\bar{1} \\ \bar{1} & \mapsto & \bar{1}1
    \end{array}
\]
on the two-letter alphabet $\{1,\bar{1}\}$; see \cite{AS,QBook} for
background. Let $v=v^{}_{0}v^{}_{1}v^{}_{2}\ldots$ be the unique
one-sided fixed point of $\varrho$ with $v^{}_{0}=1$. This infinite
word satisfies
\begin{equation}\label{eq:tmrec}
    v^{}_{2n}=v^{}_{n} \quad\text{and}\quad
    v^{}_{2n+1}=\bar{v}^{}_{n} 
\end{equation}
for all $n\in\mathbb{N}_{0}^{}$, where $\bar{\bar{1}}=1$. Consider the
bi-infinite word $w$, defined as
\[
   w_{n}= \begin{cases} v_{n} , & \text{for $n\ge 0$} , \\
                        v_{-n-1} , & \text{for $n<0$} ,
   \end{cases}
\]
which is the unique reflection-symmetric fixed point of $\varrho^{2}$
with (admissible) seed $w^{}_{-1}w^{}_{0}=11$. It is cube-free and
thus aperiodic. If $S$ denotes the (two-sided) shift, defined by
$(Sw)^{}_{n}=w^{}_{n+1}$, the Thue-Morse hull is the compact space
$\mathbb{X}=\overline{\{ S^{m}w\mid m\in\mathbb{Z}\} }$, where the
closure is taken in the obvious product topology. Note that
$\mathbb{X}$ coincides with the (discrete) local indistinguishability
(or LI) class of $w$, and that the topological dynamical system
$(\mathbb{X},S)$ is strictly ergodic; see \cite{QBook} for details.

Here, we consider the weighted Dirac comb
\begin{equation} \label{eq:tm-dc}
   \omega = \sum_{n\in\mathbb{Z}} w_{n}\,\delta_{n},
\end{equation}
where $\delta_{x}$ denotes the normalised Dirac measure on
$\mathbb{R}$, located at $x$, and $\bar{1}=-1$. In particular,
$\omega$ is a translation bounded, signed measure on $\mathbb{R}$. The
corresponding autocorrelation (or Patterson) measure $\gamma$,
obtained by a volume-averaged convolution \cite{Hof,BM} of $\omega$
with its reflected counterpart, reads
\begin{equation}\label{eq:auto}
    \gamma  \, = \, \eta\, \delta^{}_{\mathbb{Z}}
      \,:= \sum_{m\in\mathbb{Z}} \eta\ts(m)\,\delta_{m} \ts ,
\end{equation}
with the coefficients $\eta\ts(m)$ obtained as the limits
\[
   \eta\ts(m)\,  = \lim_{N\to\infty} \frac{1}{2N+1} 
   \sum_{\substack{-N\le k,\ell\le N\\ k-\ell=m}} 
   \!\! w^{}_{k} \, w^{}_{\ell}
    \, = \lim_{N\to\infty} \frac{1}{2N+1} 
   \sum_{n=-N}^{N} w_{n} \, w_{n-m} \ts .
\]
Note that these limits exist due to unique ergodicity, compare
\cite{QBook}, with $\eta\ts(0)=1$.  The autocorrelation coefficients
satisfy $\eta\ts(-m)=\eta\ts(m)$ and, due to \eqref{eq:tmrec}, the
recursions
\begin{equation}\label{eq:tmrecur}
     \eta\ts(2m) = \eta\ts(m) \quad\text{and}\quad
     \eta\ts(2m+1) = -\frac{1}{2} \bigl(\eta\ts(m)+\eta\ts(m+1)\bigr),
\end{equation}
valid for all $m\ge 0$. In particular, $\eta\ts(1)=-\frac{1}{3}$, and
all coefficients are uniquely specified. 

By construction, $\eta$ is a positive definite function on
$\mathbb{Z}$, wherefore the Herglotz-Bochner theorem
\cite[Thm.~I.7.6]{KBook} guarantees the existence of a finite 
positive measure $\nu$ on $[0,1)$ with
\begin{equation}\label{eq:etatm}
   \eta\ts(m) \, = \int_{0}^{1}\! e^{2\pi i my} \dd\nu\ts(y)\ts . 
\end{equation}
The diffraction measure $\widehat{\gamma}$ of the weighted Thue-Morse
comb $\omega$ of \eqref{eq:tm-dc} is the Fourier transform of
$\gamma$. An elementary calculation shows that
\[
    \widehat{\gamma} = \nu * \delta^{}_{\mathbb{Z}}
    \quad\text{and}\quad
    \nu = \widehat{\gamma}\ts\big|^{}_{[0,1)} \ts .
\]
In this formulation, the spectral properties of $\widehat{\gamma}$
follow immediately from those of $\nu$. In particular, $\nu$ has a 
unique decomposition
\[
    \nu = \nu^{}_{\text{pp}} + \nu^{}_{\text{sc}} + 
     \nu^{}_{\text{ac}}
\]
into its \underline{p}ure \underline{p}oint, \underline{s}ingular
\underline{c}ontinuous and \underline{a}bsolutely
\underline{c}ontinuous parts (relative to Lebesgue measure $\lambda$);
see \cite[Thms.~I.13 and I.14]{RS}.
 
It was first shown by Kakutani \cite{K} that $\nu$ is a purely
singular continuous measure, so $\nu=\nu^{}_{\text{sc}}$.  Let us
briefly adapt this to our setting.  By Wiener's criterion
\cite[Cor.~7.11]{KBook}, $\nu^{}_{\text{pp}}=0$ is equivalent to
$\lim_{N\to\infty}\frac{1}{2N+1}\varSigma(N)=0$, where
\[
    \varSigma(N)\, := \sum_{m=-N}^{N} 
    \bigl(\eta\ts(m)\bigr)^{2}\ts .
\]
When $N\ge 1$, we can use the recursion relations \eqref{eq:tmrecur}
to derive $\varSigma(4N)\le\frac{3}{2}\,\varSigma(2N)$.  With $\alpha
= \log_{2} (3/2) <1$, one then obtains the estimate $\frac{1}{N}
\varSigma(N) = \mathcal{O} (1/N^{1-\alpha})$. At this point, we
know that $\nu=\nu^{}_{\text{sc}}+\nu^{}_{\text{ac}}$, and define
the distribution function
\begin{equation}\label{eq:tm-stieltjes}
    F(x) := \nu\bigl([0,x]\bigr),
\end{equation}
which is a continuous function of bounded variation on
$[0,1]$. It satisfies $F(0)=0$ and $F(x)+F(1-x)=1$, the
latter as a consequence of the symmetry of $\nu$.

Following \cite{K} and viewing $\nu$ as a Lebesgue-Stieltjes measure
\cite[Ch.~X]{Lan93} with distribution function $F$, the two recursion
relations \eqref{eq:tmrecur} imply the identities
\begin{equation}\label{eq:tm-funeq}
   \begin{split}
   \dd F(\tfrac{x}{2}) + \!\dd F(\tfrac{x+1}{2})
   &\,= \dd F(x)\ts ,\\
   \dd F(\tfrac{x}{2}) - \!\dd F(\tfrac{x+1}{2})
   &\,=-\cos(\pi x) \dd F(x)\ts ,
   \end{split}
\end{equation}
for all $x\in[0,1]$. The left-hand sides are obtained from the
corresponding sides of the recursions by a change of variables,
followed by a split of the new integration region into two parts.  The
actual equality of the measures follows because we obtain equality of
the integrals over arbitrary trigonometric polynomials, whence the
Fourier uniqueness theorem \cite[Thm.~I.2.7]{KBook} applies. The
relations \eqref{eq:tm-funeq} must also hold for the two continuous
components of $\nu$ separately, because
$\nu^{}_{\text{sc}}\perp\nu^{}_{\text{ac}}$ by
\cite[Thm.~VII.2.4]{Lan93}.

Writing $\nu^{}_{\text{ac}}=g\lambda$ with $g\in
L^{1}\bigl([0,1]\bigr)$, the identities \eqref{eq:tm-funeq} result in
\[
   \tfrac{1}{2}\bigl( g(\tfrac{x}{2}) + g(\tfrac{x+1}{2})\bigr)
   = g(x)\quad\text{and}\quad
   \tfrac{1}{2}\bigl( g(\tfrac{x}{2}) - g(\tfrac{x+1}{2})\bigr)
   =-\cos(\pi x) \, g(x)\ts ,
\]
this time for almost all $x\in [0,1]$. Defining now
$\eta^{}_{\text{ac}}(m) = \int_{0}^{1} e^{2\pi i m x} g(x)\dd x$, we
have $\eta^{}_{\text{ac}}(-m)=\eta^{}_{\text{ac}}(m)$ and inherit a
set of recursions identical to \eqref{eq:tmrecur}, however with
$\eta^{}_{\text{ac}}(0)=0$ as a result of the Riemann-Lebesgue lemma
\cite[Thm.~I.2.8]{KBook}, and hence $\eta^{}_{\text{ac}} (m) =0$ for
all $m\in\mathbb{Z}$. By the Fourier uniqueness theorem
\cite[Thm.~I.2.7]{KBook}, this implies $g(x)=0$ almost everywhere, and
hence $\nu^{}_{\text{ac}}=0$. Since $\nu$ itself cannot be the zero
measure, we have $\nu=\nu^{}_{\text{sc}}\neq 0$, and the Thue-Morse
diffraction measure $\widehat{\gamma}$ is a purely singular
continuous, $\mathbb{Z}$-periodic, positive measure. It is the same
for all Dirac combs of members of the hull $\mathbb{X}$.

The remainder of this note is concerned with an explicit calculation
of $\nu$, via its distribution function $F$, which does not seem to
have been calculated in the literature so far.  Adding the two equations
for $\!\dd F$ from \eqref{eq:tm-funeq}, followed by integration,
results in the functional equation
\begin{equation}\label{eq:tm-inteq}
    F(x) = \frac{1}{2}\int_{0}^{2x}\bigl(1-\cos(\pi y)\bigr)\dd F(y)
\end{equation}
for the continuous function $F$, valid for $x\in [0,\frac{1}{2}]$. 
Since $F$ is a continuous non-decreasing function on $[0,1]$, the
difference $f(x)=F(x)-x$ defines a continuous function $f$ of bounded
variation, with $f(0)=0$ and $f(x)+f(1-x)=0$ for all $x\in [0,1]$. As
such, $f$ possesses a uniformly converging Fourier series of the form
\[
    f(x) = \sum_{m\ge 1} b_{m}\ts \sin(2\pi m x)\ts ,
\]
with $b_{m}=2\int_{0}^{1}\sin(2\pi m x)\,f(x)\dd x$. {}From
Eqs.~\eqref{eq:etatm} and \eqref{eq:tm-stieltjes}, one finds
$\eta\ts(m) = \pi\ts m\, b_{m}$, so that
\begin{equation}\label{eq:tmsol}
    F(x) = x + 
    \sum_{m\ge 1} \frac{\eta\ts(m)}{m\pi} \sin(2\pi mx)
\end{equation}
is a uniformly converging Fourier series representation of the
Thue-Morse distribution function. This function satisfies the symmetry
relation $F(x)+F(1-x)=1$ for all $x\in [0,1]$, and is a solution of the
integral equation \eqref{eq:tm-inteq}.

To interpret \eqref{eq:tm-inteq} in a wider setting, let us introduce
the space
\begin{equation}\label{eq:tm-space}
    D = \Bigl\{ G \in C\bigl([0,1],\mathbb{R}\bigr) 
    \;\Big|\; \substack{\text{$G(0)=0$, 
    $G$ non-decreasing and}\\[0.3mm] \text{$G(x)+G(1-x)=1$ on $[0,1]$}}
    \Bigr\}
\end{equation}
of continuous, non-decreasing distribution functions on $[0,1]$ with
the required symmetry, and define the mapping $\varPhi:\,
D\longrightarrow D$ via $G\mapsto\varPhi\ts(G)$ with
\[ 
    \bigl(\varPhi\ts(G)\bigr)(x) = 
    \begin{cases}
      \frac{1}{2}\int_{0}^{2x}\bigl(1-\cos(\pi y)\bigr) \dd G(y) \ts ,
       & \text{if $0\le x\le\frac{1}{2}$}\ts , \\
       1-\bigl(\varPhi\ts(G)\bigr)(1\!-\!x)\ts , 
       & \text{if $\frac{1}{2}<x\le 1$}\ts .  
    \end{cases}
\]
One can show that $F$ of \eqref{eq:tmsol} is the only solution of
$\varPhi\ts(F)=F$ within $D$, because $\varPhi$ is a weak contraction
with respect to the metric $d(G,H):=V(G-H)$, where $V$ denotes the
total variation of continuous functions on $[0,1]$; see
\cite[Ch.~X]{Lan93} for background.

For a numerical high-precision computation of $F$, we employ a
Volterra-type iteration within $D$, which is superior to using
\eqref{eq:tmsol} directly.  Starting from $F_{0}(x)=x$, where
$F_{0}\in D$, we set $F^{}_{n+1} = \varPhi\ts(F_{n})$ for $n\ge 0$.
This defines the (uniformly converging) sequence of distribution
functions in $D$ given by
\begin{equation}\label{eq:ftmexplicit}
   F_{n}(x) = x + \sum_{m=1}^{2^{n}-1}
   \frac{c^{(n)}_{m}}{m\pi}\, \sin(2\pi mx)\ts .
\end{equation}
Here, the coefficients $c^{(n)}_{m}$ satisfy the recursions
\[
   c^{(n+1)}_{2m}= c^{(n)}_{m}\quad\text{and}\quad
   c^{(n+1)}_{2m+1}= -\frac{1}{2}\bigl(c^{(n)}_{m}+c^{(n)}_{m+1}\bigr)
\]
for $n\ge 1$ and $0 \le m \le 2^n - 1$, together with $c^{(1)}_{0}=0$,
$c^{(1)}_{1}=-\frac{1}{2}$ and $c^{(n+1)}_{1}=
-\frac{1}{2}\bigl(1+c^{(n)}_{1}\bigr)$. They approach the
autocorrelation coefficients $\eta\ts(m)$ via $\lim_{n\to\infty}
c^{(n)}_{m} = \eta\ts(m)$ for arbitrary $m\in\mathbb{N}$.

\begin{figure}
\centerline{\includegraphics[width=0.8\textwidth]{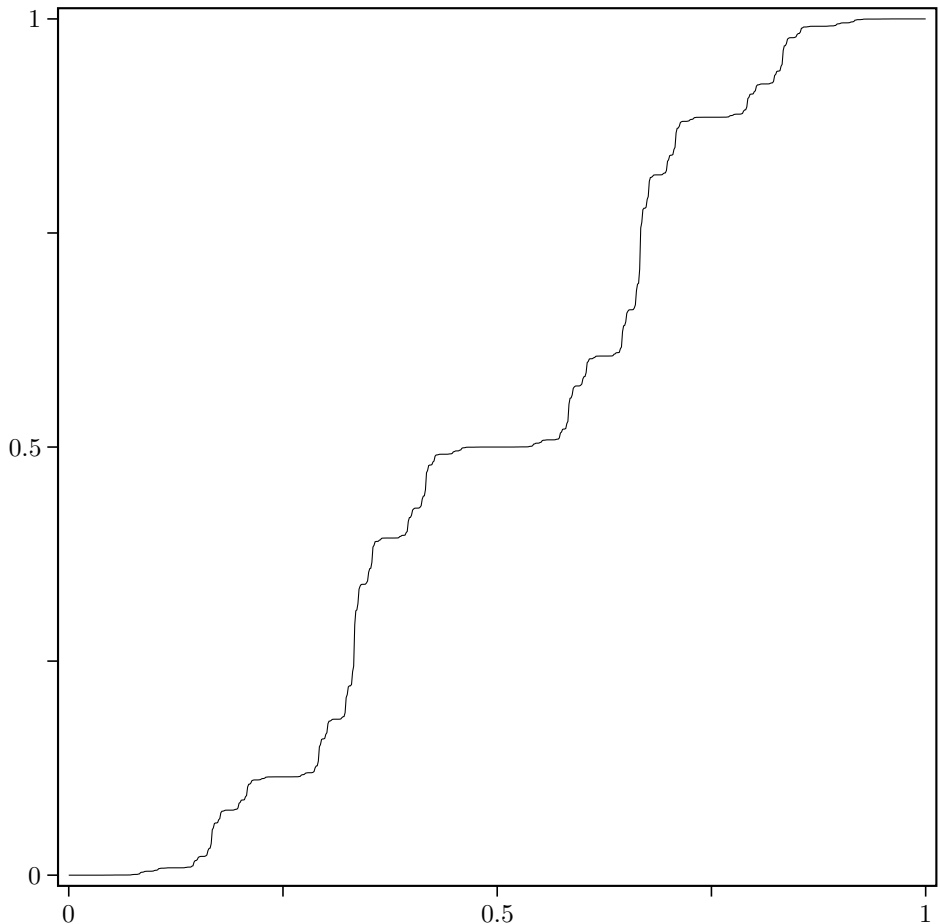}\quad}
\caption{Distribution function $F$ of the 
purely singular continuous diffraction measure of the Thue-Morse chain
in the interval $[0,1]$. Note that $F$ is a strictly increasing
function, despite the appearance of plateau-like regions.}
\label{fig:ftm}
\end{figure}

By construction, the distribution functions $F_{n}$ represent
absolutely continuous measures (although their limit does not).
Writing $\dd F_{n}(x)=f_n(x)\dd x$, one finds $f^{}_{0}(x)=1$, and the
functional equation \eqref{eq:tm-inteq} induces the recursion
\[
    f_{n+1}(x) \,=\, \bigl(1-\cos(2\pi x)\bigr) f_{n}(2x)
               \,=\, 2\ts\bigl(\sin(\pi x)\bigr)^{2} f_{n}(2x)
\]
for $n\ge 0$. This gives the well-known explicit representation 
\begin{equation}\label{eq:tm-riesz}
   f_{n}(x) \,= \prod_{\ell=0}^{n-1}
       \bigl(1-\cos(2^{\ell+1}\pi x)\bigr)
       \,=\, 2^{n}\prod_{\ell=0}^{n-1}
       \bigl(\sin(2^{\ell}\pi x)\bigr)^{2}
\end{equation}
of the Thue-Morse measure as a Riesz product; 
compare \cite[Sec.~1.4.2]{Que95}.

The Volterra iteration leads to a sequence
$(F_{n})_{n\in\mathbb{N}_{0}^{}}$ of continuous distribution functions
that converge uniformly (on $[0,1]$) to $F$, which is continuous as
well. The latter is shown in Figure~\ref{fig:ftm} and resembles the
classic middle-thirds Cantor measure in various aspect, though it has
full support (see below).  Note that the corresponding sequence
$(\!\dd F_{n})_{n\in\mathbb{N}_{0}^{}}$ of absolutely continuous
measures is only vaguely convergent, with the limit being purely
singular continuous.  Therefore, it is somewhat misleading to show a
density for the Thue-Morse measure, as is often found in the
literature. Still, it may be instructive to inspect the sequence of
densities $f_{n}$ in order to get some intuition on the singular
nature of the Thue-Morse measure, or to study some of its scaling
properties; see \cite{GVG} and references therein.

On initial inspection, the distribution function $F$ seems to have a
plateau around $x=\frac{1}{2}$, with exact value $\frac{1}{2}$. More
generally, as suggested by the Riesz products \eqref{eq:tm-riesz}, one
might expect a plateau around any $x\in\bigl(\{0,1\}\cup\{m/2^{k}\mid
k\in\mathbb{N},\, m\in\mathbb{N}\text{ odd}\}\bigr)\cap[0,1]$, because
the densities $f_{n}$ vanish at $x$ for all sufficiently large $n$,
with the order of this zero linearly increasing with $n$. However,
these potential gaps are all closed (see below).  The corresponding
values of $F$ can be calculated with the series expansion
\eqref{eq:tmsol}.  Unlike the situation of the gap labelling theorem
for one-dimensional Schr\"{o}dinger operators \cite{BHZ00}, where one
knows the values of the integrated density of states (IDOS) on the
(always non-overlapping) plateaux but not their positions, we know the
possible locations of the plateaux, but cannot see a topological
constraint for the corresponding values of the distribution function
(which can be determined from \eqref{eq:tmsol}).

It is interesting to note that the set of potential plateau locations
coincides with the set of potential (but in our case extinct) Bragg
peak positions, so the (extinct) Bragg peaks appear to `repel' the
continuous diffraction spectrum. Nevertheless, one has
\[
    \mathrm{supp} (\dd F)  =  \mathrm{supp} (\nu)
    = [0,1]
\]
by \cite[Prop.\ 28]{BBK}, which also implies that $F$ is a strictly
increasing function. In particular, as one can immediately see
from \eqref{eq:tm-inteq}, $F(x)=0$ forces $F(2x)=0$, which (when
repeated) contradicts $F(1/2)=1/2$ unless $x=0$. This shows that there
is no gap around $0$ (and none around $1$ by symmetry). The general
argument uses that $\nu$ is a regular Borel measure; see \cite{BBK}
and references therein for details.

The methods used above can also be applied to other substitutions of
constant length that fail to have pure point spectrum. This can be
decided on the basis of Dekking's criterion; see \cite[Sec.~6]{Que95}
for details. Although one still has to check Wiener's criterion and to
find the analogue of the functional equation \eqref{eq:tm-funeq}, this
approach seems worth pursuing.\bigskip

\bigskip

\noindent \textit{Acknowledgements.}  It is a pleasure to thank J.\
Bellissard, N.P.\ Frank, R.V.\ Moody and B.\ Solomyak for discussions,
and the School of Mathematics and Physics at the University of
Tasmania for their kind hospitality.  This work was supported by the
German Research Council (DFG), within the CRC 701, and by EPSRC via
Grant EP/D058465.

\end{document}